\documentclass[12pt,reqno]{amsart}

\usepackage{color}
\usepackage{amsmath, amsthm, amssymb}
\usepackage{amsfonts}
\usepackage[ansinew]{inputenc}
\usepackage[dvips]{epsfig}
\usepackage{graphicx}
\usepackage[english]{babel}
\usepackage{hyperref}

\usepackage{cite}
\usepackage{graphicx}
\usepackage{amscd}
\usepackage{color}
\usepackage{bm}
\usepackage{enumerate}

\usepackage{verbatim}
\usepackage{hyperref}
\usepackage{amstext}
\usepackage{latexsym}
\theoremstyle{plain}
\newtheorem{Theorem}{Theorem}[section]

\newtheorem{Proposition}[Theorem]{Proposition}
\newtheorem{Lemma}[Theorem]{Lemma}
\newtheorem{Remark}[Theorem]{Remark}

\numberwithin{Theorem}{section}
\numberwithin{equation}{section}

\def\proof{\noindent{{\bf Proof. }}}
\def\square{\vbox{
\hrule height .4pt \hbox{\vrule width .4pt height 7pt \kern 7pt
\vrule width .4pt} \hrule height .4pt }}
\def\QED{\hfill {$\square$}\goodbreak \medskip}

\newcommand{\average}{{\mathchoice {\kern1ex\vcenter{\hrule height.4pt
width 6pt depth0pt} \kern-9.7pt} {\kern1ex\vcenter{\hrule
height.4pt width 4.3pt depth0pt} \kern-7pt} {} {} }}

\def\R{\mathbb{R}}

\renewcommand{\a }{\alpha }

\newcommand{\e }{\varepsilon }

\renewcommand{\l }{\lambda }

\newcommand{\vp }{\varphi }

\renewcommand{\t }{\tau }

\renewcommand{\o }{\omega }
\renewcommand{\O }{\Omega }

\newcommand{\ov}{\overline}

\newcommand{\be}{\begin{equation}}
\newcommand{\ee}{\end{equation}}

\newcommand{\al}{\alpha}

\newcommand{\calL }{\mathcal{L}}

\newcommand{\calF}{{\mathcal F}}

\newcommand{\N}{\mathbb{N}}

\newcommand{\cH}{{\mathcal H}}

\newcommand{\cL}{{\mathcal L}}

\DeclareMathOperator{\id}{id}

\renewcommand{\epsilon}{\varepsilon}

\DeclareMathOperator{\spann}{span}
\begin{document}

\title[ An  overdetermined problem  for sign-changing  eigenfunctions ]
{An  overdetermined problem  for sign-changing  eigenfunctions in unbounded domains}

%\author{Mouhamed Moustapha Fall}
%\address{M. M. F.: African Institute for Mathematical Sciences in Senegal, KM 2, Route de
%Joal, B.P. 14 18. Mbour, Senegal.}
%\email{mouhamed.m.fall@aims-senegal.org}

\author{Ignace Aristide Minlend}
\address{Department of Quantitative Techniques, Faculty of Economics and Applied Management, University of Douala, BP. 2701 - Douala.}
\email{ignace.a.minlend@aims-senegal.org}
\email{ignace.minlend@univ-douala.com}

%\author{Tobias Weth}
%\address{T.W.:  Goethe-Universit\"{a}t Frankfurt, Institut f\"{u}r Mathematik.
%Robert-Mayer-Str. 10 D-60629 Frankfurt, Germany.}

%\email{weth@math.uni-frankfurt.de}

\keywords{Overdetermined problems,  Sign-changing, Bifurcation}

%\subjclass[2010]{Primary ; Secondary}

\begin{abstract}
We study  the existence of non-trivial unbounded domains of $\O\subset \R^2$  where the equation
\begin{align}
  - \lambda u_{xx} -u_{tt} &=  u  \qquad \text{in $\Omega$,}\nonumber\\
   u &=0 \qquad \text{on $\partial \Omega$,}\nonumber
\end{align}
is solvable subject to the conditions
\begin{align}
 \frac{\partial u}{\partial \eta}  =-1\quad \text{on $\partial \Omega^+$} \quad \textrm{and}\quad  \frac{\partial u}{\partial \eta}  =+1\quad \text{on $\partial \Omega^-$.}\nonumber
\end{align}
For every integer $m\geq 0$, we prove the existence of a  family of unbounded domains  $\O\subset \R^2$  indexed by $0 \leqslant\ell\leqslant 2m$, where the above problem admits  periodic  sign-changing  solutions.  The domains we construct are periodic in the first  coordinate in $ \R^2$, and they bifurcate from  suitable  strips. 
\end{abstract}
\maketitle

\textbf{MSC 2010}:  35J57, 35J66, 35N25, 58J55, 35J25 
\maketitle
\section{Introduction and main result}

This paper  is concerned with  the  existence of periodic  sign-changing solutions to  a particular  prototype of overdetermined boundary value problem on strip domains of  the plane. In  1971, Serrin \cite{Serrin}  proved   by   Alexandrov \cite{Alexandrov} moving plane method that the only bounded and
regular domains  in the Euclidean space  $\mathbb{R}^N$  where the overdetermined  problem  \begin{equation}\label{eq:ProA3}
 -\Delta u=1\quad \text{ \quad in $\Omega$}
\end{equation}
and 
\begin{equation}\label{eq:ro3}
u=0, \quad \partial_\nu u=\textrm{const} \quad \text{ on \quad $\partial \Omega$}
\end{equation}
is solvable are balls. Here $\nu$ is the unit outer to the boundary. \\
Soon after  this  celebrate result  was  communicated to the PDE community, several authors have developed interest in  the study of symmetry properties   as well as rigidity  results of the  overdetermined problem \eqref{eq:ProA3}-\eqref{eq:ro3} including  the more general equation
 \begin{align}\label{probcafa}
-\Delta u = f(u) \quad \text{in $\Omega$}, \qquad u=0, \quad \partial_\nu u=\textrm{const}  \qquad \text{on $\partial \Omega$},
\end{align} where  $f: [0,\infty) \to \R$ is  a locally Lipschitz function.  We refer the reader to    
\cite{Alessandrini, Gazzola, Greco,  Philippin,  GarofaloLewis, Prajapat, PaynePhilippin, Lamboley, FragalaGazzola, BrockHenrot, Reichel, PhilippinPayne, BerchioGazzolaWeth, FraGazzolaKawohl}, \cite{farina-valdinoci,farina-valdinoci:2010-1,farina-valdinoci:2010-2,farina-valdinoci:2013-1, farina-valdinoci:2013-2, BCNI}.  In the classical literature a smooth domain $\O$ where the problem \eqref{probcafa} is solvable is called \emph{$f$-extremal domain.}. In 1997,  Berestycki, Caffarelli and Nirenberg \cite{BCNI} addressed the classification of  unbounded \emph{$f$-extremal domains} and conjectured that if  $\Omega \subset \R^N$ is an \emph{$f$-extremal domain} such that  $\R^N \setminus \overline \Omega$  is
connected   and  \eqref{probcafa} admits a bounded solution, then $\Omega$ is either a half plane,  or a generalized  cylinder $\R^j \times B$  where $B$ is the unit Euclidean ball in
 $\R^{N-j}$, or  the complement $B^c$  of a ball $B \subset
\R^N$  or a cylinder. 
For  $f(u)=\lambda_1 u$,  where  $\lambda_1$
is the first eigenvalue of the Laplacian with $0$-Dirichlet boundary
condition, this conjecture was disproved in dimension  $N\geq 3$ by
Sicbaldi \cite{Sic},  and later in dimension $N\geq2$ by  Sicbaldi
and  Schlenk in \cite{ScSi}, where they proved existence of periodic
and unbounded extremal domains  bifurcating from straight cylinder
$\R \times B$.  Subsequently,  we   considered the case  $f\equiv 1$ in  \cite{ Fall-MinlendI-Weth} and  proved the existence of periodic unbounded domains bifurcating from generalized-type cylinder domains in $\R^N$. We also refer  to   \cite{Fall-MinlendI-Weth2, Ku-Pra, morabito-sicbaldi}, where  the  conjecture in \cite{BCNI}  is  addressed in space forms. \\

We note that  the results  in  the previous  works  all assume  one-sign solutions, while the existence of sign-changing solutions in the context of overdetermined boundary value problems is a subject of intensive investigation in the literature. 

In the present paper, we deal with  the existence of sign-changing solutions  to  a specific prototype of overdetermined boundary value problem given by
\begin{equation} \label{eq:perturbed-strip2}
  \left \{
    \begin{aligned}
        - \lambda u_{xx} -u_{tt} &=  u && \qquad \text{in $\Omega$,}\\
              u &=0 &&\qquad \text{on $\partial \Omega$,}\\
             \frac{\partial u}{\partial \eta}  &=-1 &&\qquad \text{on $\partial \Omega^+$,}\\
              \frac{\partial u}{\partial \eta}  &=+1 &&\qquad \text{on $\partial \Omega^-$}
    \end{aligned}
       \right.
     \end{equation}
for  some $\lambda >0$. 
Here, $\nu$ is the outer unit normal vector field to the boundary of  $\Omega$,$$ \partial \Omega^+= \{ (x, t) \in \partial \Omega, \quad t> 0  \}\quad \textrm{and} \quad\partial \Omega^-= \{ (x, t) \in \partial \Omega,  \quad t< 0  \}.$$ 
Our aim is to  prove the existence of sign-changing solutions to  problem  \eqref{eq:perturbed-strip2} on unbounded domains $\O \subset \R^2$. 

Related to this paper are  the following  conjectures of Schiffer in spectral  theory, see \cite{schiffer:1954, schiffer:1955, ChatelainHe}. 

\textbf{\emph{Conjecture 1:}}
Let $\O$ be a simply connected  domain in $\R^N$ and $\lambda\ne 0$.  Then there exists a  Dirichlet eigenfunction $u\ne 0$ satisfying 
\begin{equation}
(P_{\lambda}):  \label{eq:schiffer}
  \left \{
    \begin{aligned}
        -\Delta u &= \lambda u && \qquad \text{in $\Omega$,}\\
           u &=0&& \qquad \text{on $\partial \Omega$,}\\
           \frac{\partial u}{\partial \eta}  &=const. &&\qquad \text{on $\partial \Omega$.}
    \end{aligned}
       \right.
\end{equation}
if and only if $\O$ is a ball.

\textbf{\emph{Conjecture 2:}}
Let $\O$ be a simply connected  domain in $\R^N$. If  $\lambda >0$, then  there exists  $u\ne 0$ such that  
\begin{equation}
(Q_{\lambda}):  \label{eq:schiffer}
  \left \{
    \begin{aligned}
        -\Delta u &= \lambda u && \qquad \text{in $\Omega$,}\\
              \frac{\partial u}{\partial \eta} &=0&&\qquad \text{on $\partial \Omega$,}\\
           u  &=const. &&\qquad \text{on $\partial \Omega$.}
    \end{aligned}
       \right.
\end{equation}
if and only if $\O$ is a ball.

Indeed if $u$ is a solution of \eqref{eq:perturbed-strip2}, then the function 
$$
W^{\lambda} (y, \zeta):= u(\lambda y, \sqrt{\lambda}\zeta)
$$
solves  the  overdetermined problem
\begin{equation} \label{eq:W1}
  \left \{
    \begin{aligned}
        - \Delta W^{\lambda}  &= \lambda W^{\lambda} && \qquad \text{in $\Omega_\lambda$,}\\
              W^{\lambda}   &=0 &&\qquad \text{on $\partial \Omega_{\lambda}  $,}
    \end{aligned}
       \right.
     \end{equation}  
     \begin{align}\label{eqNeuW1}
       \frac{\partial W^{\lambda}}{\partial \eta}  &=-1 \quad \text{on  \quad $\partial  \Omega_\lambda^+$} \quad \textrm{and}\quad \frac{\partial W^{\lambda}}{\partial \eta} =+1 \qquad \text{on \quad  $\partial  \Omega_\lambda^+$,}
\end{align}
where $$
\Omega_\lambda:= \{(\frac{x}{\lambda}, \frac{\t}{\sqrt{\lambda}}),\: \quad \: (x,\t) \in  \Omega\}.
$$
The solution  $W^{\lambda}$  does not assume a  constant Neumann data  at  boundary $\Omega_\lambda$ as  \eqref{eqNeuW1} reveals.  Nevertheless the problem  \eqref{eq:W1}-\eqref{eqNeuW1} can be seen as a refined  version of  conjecture 1 above. In contract, interchanging the boundary conditions in \eqref{eq:perturbed-strip2} and  considering  constant zero Neumann  data and constant Dirichlet  equals  $-1$,  we are led to 
the problem
\begin{align}\label{eq:Wl4}
  \left \{
    \begin{aligned}
        - \Delta W^{\lambda}  &= \lambda W^{\lambda} && \qquad \text{in $\Omega_\lambda$,}\\
              W^{\lambda}   &=-1 &&\qquad \text{on $\partial \Omega_{\lambda}  $,}
    \end{aligned}
       \right.
     \end{align}  
     \begin{align}\label{eqNeuWl4}
       \frac{\partial W^{\lambda}}{\partial \eta}  &=0\quad \text{on  \quad $\partial  \Omega_\lambda$,}\end{align}
which is a prototype in the  conjecture 2.  We announce here our forthcoming paper where we are  addressing \eqref{eq:Wl4}-\eqref{eqNeuWl4}. We alert that this  problem leads to a serious  loss of regularity which we hope  to overcome soon.\\
       
Up to now, only  few works have been carried out regarding the existence  on sign-changing solutions of overdetermined boundary value problems and the most recent  in this context are \cite{Ramm, Deng, BCanutoDRial, B.Canuto}.  In   \cite{B.Canuto}  and \cite{BCanutoDRial}, the authors  addressed   the following  question: for $\o\in \R$, is it true that the only domain such that there exists a solution to the overdetermined problem
\begin{equation}\label{eq:SignchangeBounded}
 \begin{cases}
    \Delta u+\o^2 u=-1 & \quad \textrm{in} \quad  \Omega \vspace{3mm}\\
u=0 &  \quad\textrm{on }\quad \partial\Omega
  \end{cases}
\end{equation}
and 
\begin{equation}\label{eq: SignchangeBoun2} 
 \quad   \partial_{\nu} u= \textrm{c } \quad \textrm{ on } \quad \partial \Omega.
\end{equation}
is a ball ?

Moreover, the authors  in  \cite{BCanutoDRial}  proved  the  under suitable assumptions on $\o$ that the  only bounded domain $\O$ such that there exists a solution  to \eqref{eq:SignchangeBounded}-\eqref{eq: SignchangeBoun2} is the ball $B_1$, independent  on the sign of $u$, provided  $\partial \O$ is a perturbation of the unit sphere $\partial B_1$ in $\R^N$. A similar result   was derived  in \cite{B.Canuto} by considering a different Neumann boundary condition. \\

Concerning the case of unbounded domains, our  paper  appears to be the first and in fact, we  are not aware of any other existing  result for sign-changing solutions in unbounded domains arising in the context of  overdetermined  boundary value problems. \\

To state our main result, we  consider the reference domain $$\Omega_*= \R \times (-(2 m+1) \pi, (2 m+1) \pi) \subset \R^2$$ for some fixed non-negative integer $m$.  Then  problem \ref{eq:perturbed-strip2} is solved  on  $\Omega_*$ 
by  the  ($\lambda$-independent) solution  $u_*(x,t)= \sin t$ for all $\lambda>0$. Our aim is construct domains $\Omega$ close to $\Omega_*$ with the property that the overdetermined  boundary value problem \ref{eq:perturbed-strip2}  is solvable on $\Omega$ by some sign-changing solution $u$.

We consider the open set $$
Y_2^+:= \{h \in C^{2,\alpha}_{p, e}(\R)\::\: h> 0 \quad \text{$h$ is even in $x$}\},$$
where  $C^{2,\alpha}_{p, e}(\R)$ stands for the space of  $2\pi$ periodic functions in $\R$. For  a function $h \in Y_2^+ $, we  define  the domain               
\begin{equation}\label{eq: domainh}
\Omega_h:= \{(x, \frac{\t}{h(x)})\::\: (x,\t) \in  \Omega_*\}.
\end{equation}
Our main result states the existence of  sign-changing solutions  $u$  to  the overdetermined problem  \eqref{eq:perturbed-strip2} on domains of the form $\Omega = \Omega_h$.

\begin{Theorem}\label{Strip1}
Let $m\in \mathbb{N}$ and $\alpha \in (0,1)$. There exists a strictly decreasing  and finite  sequence  $(\mu_\ell(m))_{0\leq \ell \leq 2m}$  of real numbers in  $(0, 1)$ with the following properties: for every ${0\leq \ell \leq 2 m}$, there  exists  ${\e_\ell}>0$ and a smooth curve
$$
(-{\e_\ell},{\e_\ell}) \to   (0,+\infty) \times C^{2,\alpha}(\overline {\Omega_{*}}) \times C^{2,\alpha}(\R),\qquad s \mapsto (\lambda_\ell(s),\varphi^\ell_s, \psi^\ell_s)
$$
with
$ \lambda_\ell(0)= \mu_\ell(m)$, $ (\varphi^\ell_0, \psi^\ell_0)\equiv(0,0)$ and such that for every $s\in (-{\e_\ell},{\e_\ell})$, there 
exists a solution $u^\ell_s \in C^{2,\alpha}(\overline {\Omega_{ 1+\psi^\ell_s}})$ of  the overdetermined problem
\begin{equation} \label{eq:solved}
  \left \{
    \begin{aligned}
        - \lambda^\ell_s u_{xx} -u_{tt} &=  u && \qquad \text{in $ \Omega_{ 1+\psi^\ell_s} $,}\\
              u &=0 &&\qquad \text{on $\partial \Omega_{ 1+\psi^\ell_s}$,}\\
             u_\nu &=-1 &&\qquad \text{on $\partial \Omega_{ 1+\psi^\ell_s}^+$,}\\
              u_\nu &=+1 &&\qquad \text{on $\partial  \Omega_{ 1+\psi^\ell_s}^-$}.
    \end{aligned}
       \right.
     \end{equation}
Moreover, the function  $\psi^\ell_s\in C^{2,\alpha}(\R)$ is even and $2\pi$ periodic and 
\begin{align}\label{eqa2}
\psi^\ell_s= s \bigl( \cos (x) + \kappa^\ell_s\bigr).
\end{align}
Furthermore, setting 
$$ \tilde u^\ell_s := \sin(t)+s \bigl(v_\ell + \mu^\ell_s\bigr),$$
where $$v_\ell (x,t)=  \Bigl((-1)^\ell(2m+1)\pi \sin (\sqrt{1-\mu_\ell(m)}t)+t\cos (t)\Bigl)\cos (x),$$
the solution    $u^\ell_s \in C^{2,\alpha}(\overline {\Omega_{ 1+\psi^\ell_s}})$ to  \eqref{eq:solved} is of  the form
\begin{align}\label{eqasolutionform}
u^\ell_s(x,t)= \tilde u^\ell_s\Bigl(x ,(1+ \psi^\ell_s(x))t\Bigl),
\end{align}
with a smooth curve
$$
(-{\e_0},{\e_0}) \to C^{2,\alpha}(\overline {\Omega_{*}}) \times C^{2,\alpha}(\R), \qquad s \mapsto  (\mu^\ell_s, \kappa^\ell_s)
$$ satisfying 
\begin{align*}
\int_{(-\pi,\pi)\times(-(2 m+1) \pi,(2 m+1) \pi) } \mu^\ell_s v_\ell(x,t)\,dxdt+\int_{(-\pi,\pi)}  \kappa^\ell_s \cos (x)dx=0.\\
\end{align*}
\end{Theorem}

The fact that the solution $u^\ell_s \in C^{2,\alpha}(\overline {\Omega_{ 1+\psi^\ell_s}})$ to  \eqref{eq:solved} changes sign is precise in Section \eqref{eq.setting} below. Indeed $\Omega_{ 1+\psi^\ell_s}$ is symmetric in the $t$  coordinate  and $u^\ell_s $ is odd with respect to this variable. Furthermore, the  sequence  $(\mu_\ell(m))_{0\leq \ell \leq 2m}$  whose existence is claimed in Theorem \ref{Strip1} is explicitly given by
\begin{align*}
 \mu_\ell(m)&=1-\frac{1}{4}\Bigl(   \frac{1+2\ell }{1+2m}\Bigl)^2.
\end{align*}
In particular $ \mu_0(0)= \mu_m(m)=\frac{3}{4}.$\\

We now describe the proof of Theorem \ref{Strip1}  while presenting the contents  of the paper.\\

The proof of Theorem \ref{Strip1}  is achieved by the use of Crandall-Rabinowitz bifurcation theorem, \cite{M.CR}. Our aim is  solve  the problem \eqref{eq:perturbed-strip2}  on the domain  $\O_h$  given by \eqref{eq: domainh}. In Section \ref{eqtranf}, we transform \eqref{eq:perturbed-strip2}  to the fixed domain $\O_*$, see \eqref{eq:Proe1-ss1}. Next we write  \eqref{eq:Proe1-ss1}  into the solvability  of  a  bifurcation equation  $F_\lambda(u,h)=0$  defined in \eqref{eq to bifurcate}, and   we  give the expression of  linearized operator $D F_\lambda(0,0)$ in  \eqref{eq diffope2}. In Section  \ref{exkus}, we compute the kernel of  $D F_\lambda(0,0)$ and derive the spectral properties of the operator $D F_\lambda(0,0)$ in Lemma \ref{lemorthoKer} and Proposition \ref{propPhlinearop} as well as all the preliminary  assumptions of  Crandall-Rabinowitz bifurcation theorem\cite{M.CR}. The proof of  Theorem \ref{Strip1} is completed in Section \ref{eq:ProofTheo1}. In section \ref{eq: Cradal Rabi1}, we give the  Crandall-Rabinowitz bifurcation theorem for the reader convenience. 

\bigskip
\noindent \textbf{Acknowledgements}: 
This work was  carried  out when  the author  was visiting the Institute of Mathematics and Informatics of the Goethe University Frankfurt as a Humboldt postdoctoral fellow. He is gratefully to the Humboldt Foundation for funding his research  and wishes to thank Department of Mathematics  of  the Goethe-University Frankfurt  for the hospitality. The author  also thanks his host Prof. Tobias Weth and Prof. Mouhamed Moustapha Fall for their  helpful suggestions  and  comments  throughout
the writing of this paper.

\section{The tranformed problem and its linearization}\label{eqtranf}
For some fixed non-negative integer $m$,  we consider the domain $$\Omega_*= \R \times (-(2 m+1) \pi, (2 m+1) \pi) \subset \R^2$$  and  the open set $$
Y_2^+:= \{h \in C^{2,\alpha}_{p, e}(\R)\::\: h> 0 \quad \text{$h$ is even in $x$}\},$$
where  $C^{2,\alpha}_{p, e}(\R)$ stands for the space of  $2\pi$ periodic functions in $\R$.

For  a function $h \in Y_2^+ $, we  define  the domain               
\begin{equation}\label{eq: domainh2}
\Omega_h:= \{(x, \frac{\t}{h(x)})\::\: (x,\t) \in  \Omega_*\}.
\end{equation}

We are looking for a solution $u$ of the overdetermined problem 
\begin{equation}\label{eq:perturbed-strip4}
  \left \{
    \begin{aligned}
        - \lambda u_{xx} -u_{tt} &=  u && \qquad \text{in $\Omega_h$,}\\
              u &=0 &&\qquad \text{on $\partial \Omega_h$,}\\
             \frac{\partial u}{\partial \eta}  &=-1 &&\qquad \text{on $\partial \Omega_h^+$,}\\
              \frac{\partial u}{\partial \eta}  &=+1 &&\qquad \text{on $\partial \Omega_h^-$}
    \end{aligned}
       \right.
     \end{equation}
for  some $\lambda >0$. Here, $\nu$ is the outer unit normal vector field to the boundary of  $\Omega_h$,$$ \partial \Omega_h^{+}= \{ (x, t) \in \partial \Omega_h, \quad t> 0  \}\quad \textrm{and} \quad\partial \Omega_h^{-}= \{ (x, t) \in \partial \Omega_h,  \quad t< 0  \}.$$ 
Observe that  $\Omega_h$ is parametrized by the mapping
$$ \Psi_h: \Omega_*  \to  \Omega_h , \quad  (x, \t)  \mapsto (x, t)=(x, \frac{\t}{h(x)}),$$ with inverse given by 
$ \Psi^{-1}_h:  \Omega_h   \to  \Omega_*, \quad  (x, t)  \mapsto   (x, h(x) t)$.\\

The aim is to pull back problem \eqref{eq:perturbed-strip4}  on  the fixed unperturbed domain $\Omega_*$ 
via the ansatz 
\begin{align}\label{eqans}
v(x,t)= u(x,h(x)t)=u(x,\tau) \qquad \text{for some function $u: \Omega_* \to \R$.}
\end{align}
For this we define
$$
L_\lambda:= 
\partial_{tt} + \lambda \partial_{xx}  +  \id,
$$
and we  first to compute the differential operator $L^h_\lambda$ with the property that 
\begin{align}\label{eqreladiffopets}
[L^h_\lambda u] (x,h(x)t) = [L_\lambda v](x,t) \qquad \text{for $(x,t) \in \Omega_h$}
\end{align}
for the function $v: \Omega_h \to \R$, $v(x,t)=u(x,h(x)t)$.

Setting $\xi: = (x,h(x)t)$ in the following, we compute that
$$
v_x(x,t )= u_x(\xi) + h'(x)t u_\tau (\xi),
$$
$$
v_{xx}(x,t) = u_{xx}(\xi) + 2 h'(x)t u_{x \tau }(\xi) + h''(x)t u_\tau (\xi)+t^2 [h'(x)]^2u_{\tau \tau}(\xi)
$$
and
$$
v_{tt}(x,t)= h(x)^2 u_{\tau \tau}(\xi).
$$
Consequently,
\begin{align*}
  &L_\lambda v(r,t)=  v(r,t) + \lambda v_{xx}v(r,t) + v_{tt}v(r,t)\\
  &=  u(\xi) +  \lambda u_{xx}(\xi) + \bigl(h^2(x) + \lambda t^2 [h'(x)]^2\bigr) u_{\tau \tau}(\xi) + \lambda h''(x)t u_\tau (\xi) + 2 \lambda h'(x)t u_{\tau x}(\xi)
\end{align*}
and 
\begin{align*}
L_\lambda v (x,\frac{\t}{h(x)}) &=  u(x,\tau ) +  \lambda u_{xx} (x,\tau ) + \Bigl(h^2(x) + \lambda   \t^2  \frac{h'(x)^2}{h(x)^2}  \Bigr) u_{\tau \tau} (x,\tau ) \\
  &+ \lambda \frac{ h''(x) }{h(x)} \t u_\tau (x,\tau ) + 2 \lambda  \frac{ h'(x) }{h(x)}   \tau u_{\tau x}(x,\tau ).
\end{align*}
Therefore the differential operator  $L^h_\lambda$ is given in coordinates  $(x,\tau )$  by  
\begin{align} \label{eq:reldiffope}
L_\lambda^h  = \id + \lambda \partial_{xx} +  \Bigl(h^2(x) + \lambda   \t^2  \frac{h'(x)^2}{h(x)^2}  \Bigr)  \partial_{\tau \tau}  + 2 \lambda  \frac{ h'(x) }{h(x)}   \tau \partial_{ x \tau } +  \lambda \frac{ h''(x) }{h(x)} \t \partial_\t. 
\end{align}

Next we  express the normal derivative of $u$ with respect to the outer normal vector field  on $\partial \Omega_*$  induced by the parametrization  $ \Psi_h: \Omega_*  \to  \Omega_h , \quad  (x, \t)  \mapsto   (x, \frac{\t}{h(x)})$.\\

Let the metric $g_h$ be
 defined as the pull back of the euclidean metric $g_{eucl}$ under the map $\Psi_h$, so
 that $\Psi_{h}:(\overline \O,g_h) \to (\overline \O_{h},g_{eucl})$ is an isometry.
Denote  by
$$
\eta_h: \partial  \Omega_* \to  \R^{2}
$$ the unit outer normal vector field on $\partial  \Omega_* $ with respect
to $g_h $.
Since \mbox{$\Psi_{h}: (\overline \Omega_*, g_h) \to
(\overline \O_h, g_{eucl})$} is an isometry, we have
\begin{equation}
  \label{eq:rel-mu-phi-nu-phi}
\eta_{h} = [d \Psi_h]^{-1} \mu_h \circ \Psi_h  \qquad
\text{on $\partial  \Omega_* $,}
\end{equation}
where $\mu_h: \partial \O_h \to \R^{2}$ denotes the outer
normal on $\partial \Omega_h$ with respect to the Euclidean
metric $g_{eucl}$ given by
\begin{equation}
  \label{eq:def-mu-phi}
\mu_h (x,t) = \dfrac{ 1}{\sqrt{1+ \frac{((2 m+1) \pi)^2 h'^2(x) }{ h^4(x)}}} \Bigl( \dfrac{  (2 m+1) \pi h'(x) }{ h^2(x)}, \dfrac{t}{|t|} \Bigl)  \in \R^{2} \qquad \text{for
$(x,t) \in \partial \Omega_h.$}
\end{equation}

Moreover, by (\ref{eq:rel-mu-phi-nu-phi}) we have
$\mu_h(\Psi_h (x,\tau) )=d\Psi_h (x,\tau) \eta_h (x,\tau) $ and therefore
\begin{align*}
\partial_{\eta_h} u (x,\tau) = d u (x,\tau) \eta_{h} (x,\tau) =d v (\Psi_{h}(x,\tau))  &d \Psi_{\phi}(x,\tau)  \eta_{h}(x,\tau)= d v (\Psi_{h}(x,\tau))   \mu_h(\Psi_h (x,\tau) )\\
  &= \langle  \mu_h(\Psi_h (x,\tau) ), \nabla_{ (x, t)}v (\Psi_{h}(x,\tau)) \rangle_{g_{eucl}}.
\end{align*}

From \eqref{eqans}
$$  \nabla_{ (x, t)}v (x, t)= \Bigl(u_x (x,h(x) t) + h'(x) t u_ \tau (x,h(x) t), h(x)u_\tau (x,h(x) t) \Bigl)$$
and we have 
$$  \nabla_{ (x, t)}v (\Psi_{h}(x,\tau))=  \nabla_{ (x, t)}v (x, \frac{\tau }{h(x)})=\Bigl(u_x (x,\tau) +  \frac{h'(x) }{h(x) } \tau u_ \tau (x,\tau), h(x)u_\tau (x,\tau) \Bigl)$$
and hence,
\begin{align} \label{eq:normor}
\partial_{\eta_h} u (x,\tau) = \dfrac{ 1}{\sqrt{1+ \frac{((2 m+1) \pi)^2 h'^2(x) }{ h^4(x)}}} \Bigl[ \frac{(2 m+1) \pi  h'(x) }{ h^2(x)}\Bigl(u_x (x,\tau) +  \frac{h'(x) }{h(x) } \tau u_ \tau (x,\tau) \Bigl)+  \dfrac{\tau}{|\tau|}h(x)u_\tau (x,\tau) \Bigl]
\end{align}

From \eqref{eqans} and \eqref{eqreladiffopets} the original problem  \eqref{eq:perturbed-strip4} is  equivalent to 
\begin{align}\label{eq:Proe1-ss1}
  \begin{cases}
    [L^h_\lambda u] (x, \t)= 0 &\quad \textrm{ in}\quad  \Omega_* \\
     u=0    &\quad\textrm{ on} \quad \partial  \Omega_* \\
\partial_{\eta_{h}} u \equiv -1 &\quad\textrm{ on }\quad \partial  \Omega^+_*\\
\partial_{\eta_{h}} u \equiv +1 &\quad\textrm{ on }\quad \partial  \Omega^-_*,
  \end{cases}
  \end{align}
where  $ \partial_{\eta_{h}} u$ is expressed in \eqref{eq:normor}.\\

Our aim is then to find $(u, h)$ such that \eqref{eq:Proe1-ss1} holds.
Before  proving  the existence of such a couple $(u, h)$, we make the following  which allow to reduce the  Neumann conditions in  \eqref{eq:Proe1-ss1} to a single equation. \\

\begin{Remark}\label{Rmm}
Let $u=u(x,t)$  satisfying $$ u(x,\pm(2 m+1) \pi)=Const.$$  Then  $u_x(x,\pm(2 m+1) \pi)=0$ so that  \eqref{eq:normor} yields 
\begin{align} \label{eq:nortt}
\partial_{\eta_h} u (x, +(2 m+1) \pi)& = \dfrac{ 1}{\sqrt{1+ \frac{(2 m+1) \pi)^2  h'^2(x) }{ h^4(x)}}} \Bigl[ \frac{ ((2 m+1) \pi)^2  h'^2(x) }{ h^3(x) } +  h(x) \Bigl] u_\tau (x, (2 m+1) \pi)\nonumber \\
\partial_{\eta_h} u (x,-(2 m+1) \pi) &= -\dfrac{ 1}{\sqrt{1+ \frac{((2 m+1) \pi)^2 h'^2(x) }{ h^4(x)}}} \Bigl[ \frac{ (2 m+1) \pi)^2  h'^2(x) }{ h^3(x) } +  h(x) \Bigl] u_\tau (x, -(2 m+1) \pi).\nonumber 
\end{align}
In addition to  $ u(x,\pm(2 m+1) \pi)=Const.$, if we  assume that $u=u(x,t)$  is odd in $t$, then 
$u_t (x,-(2 m+1) \pi)=u_t (x,(2 m+1) \pi)$ and 
\begin{equation}\label{eq:normlllls}
\partial_{\eta_h} u (x, +(2 m+1) \pi)=-\partial_{\eta_h} u (x,-(2 m+1) \pi).
\end{equation}

In this case  the equation 
\begin{equation}
\label{eq:Proe1-1-neumann55}
\partial_{\eta_{h}} u \equiv c \ne 0  \quad\textrm{ on }\quad \partial  \Omega_*
\end{equation}
has no solution.
Nevertheless,  one can consider the overdetermined problem
\begin{equation} \label{eq:schiffer-perturbed-for sine}
 \left \{
    \begin{aligned}
        - \lambda u_{xx} -u_{tt} &=  u && \qquad \text{in $\Omega$,}\\
           u &=0 &&\qquad \text{on $\partial \Omega$,}\\
             u_\nu &=-1 &&\qquad \text{on $\partial \Omega^+$,}\\
              u_\nu &=+1 &&\qquad \text{on $\partial \Omega^-$,}
    \end{aligned}
       \right.
     \end{equation}
where $  \partial \Omega^+= \{ (x, t) \in \partial \Omega, \quad t> 0  \}$  and  $  \partial \Omega^-= \{ (x, t) \in \partial \Omega,  \quad t< 0  \}$. It  is obvious that the  function  $\sin (t)$  solves  \eqref{eq:schiffer-perturbed-for sine} on  $\Omega_*$.

From  \eqref{eq:normlllls}, the  Neumann boundary conditions in \eqref{eq:schiffer-perturbed-for sine} therefore reduce to the single equation
\begin{align} \label{eq:normale}
\dfrac{ 1}{\sqrt{1+ \frac{((2 m+1) \pi)^2  h'^2(x) }{ h^4(x)}}} \Bigl[ \frac{ ((2 m+1) \pi)^2  h'^2(x) }{ h^3(x) } +  h(x) \Bigl] u_\tau (x, (2 m+1) \pi)=-1.
\end{align}
\end{Remark}

\section{The setting of the pull back problem  and computation of the linearized operator}\label{eq.setting}

To set up a framework for  problem   \eqref{eq:Proe1-ss1}, we  define the function spaces
$$
C^{2,\alpha}_{p,e}(\overline \Omega_*):= \{u \in C^{2,\alpha}(\overline \Omega_*)\::\: \text{$u=u (x, t)$ is  $2\pi$ periodic $x$}\},
$$
$$
X_2:= \{u \in C^{2,\alpha}_{p,e}(\overline \Omega_*)\::\: \text{$u$ is odd in $t$ and even in $x$, and $u \equiv 0$ on $\partial \Omega_*$.}\}
$$
as well as 
$$
X_0:= \{u \in C^{0,\alpha}(\overline \Omega_*)\::\: \text{$u$ is odd in  $t$ and even in $x$, and $2\pi$ periodic in $x$}\}
$$
and
$$
Y_2:= \{h \in C^{2,\alpha}_{p,e}(\R)\::\: \text{$h$ is even  in $x$}\} \quad 
Z_2:= \{z \in C^{1,\alpha}_{p,e}(\R)\::\: \text{$z$ is even  in $x$}\}.$$
We also recall  $Y_2^+:= \{h \in Y_2\::\: h>0\}$ and $u_*(x,t)= \sin t$.\\

Our aim is to prove that for  some parameter  $\lambda$, we can find the functions  $(u,h)\in  X_2 \times Y_2$  such that 

\begin{align}\label{eq:Proe11}
\left \{
    \begin{aligned}
    &[L^h_\lambda u] (x, \t)= 0 & \quad \textrm{ in}\quad  \Omega_*\\
&u=0&  \quad\textrm{ on} \quad \partial  \Omega_* \\
&\dfrac{ 1}{\sqrt{1+ \frac{(2 m+1) \pi)^2  h'^2(x) }{ h^4(x)}}} \Bigl[ \frac{ ((2 m+1) \pi)^2  h'^2(x) }{ h^3(x) } +  h(x) \Bigl] u_\tau (x, (2 m+1) \pi)=-1 & \quad \textrm{ on} \quad  \R.
  \end{aligned}
       \right.
\end{align}
Define the mapping
$$
Q: X_2 \times Y_2^+ \to  Z_2
$$ by 
$$Q (u,h):=\dfrac{ 1}{\sqrt{1+ \frac{((2 m+1) \pi)^2  h'^2(x) }{ h^4(x)}}} \Bigl[ \frac{ ((2 m+1) \pi)^2  h'^2(x) }{ h^3(x) } +  h(x) \Bigl] u_\tau (x, (2 m+1) \pi)+1.$$
and 
\begin{align}\label{eqQtilde}
\widetilde{ Q } (u,h):= Q (u+ u_*,1+h).
\end{align}
Then $\widetilde{ Q } (0,0)=0$ and by  computation
\begin{align}\label{eq:derQ}
D \widetilde{ Q }(0,0)(v,g)= v_\tau (\cdot, (2 m+1) \pi)-g(\cdot).
\end{align}
Next, we define 
\begin{align}\label{eq:maptFF}
F_\lambda: X_2 \times Y_2^+ \to X_0  \times Z_2  , \qquad F_\lambda(u,h):= (L_\lambda^{1+h} (u+u_*), \widetilde{ Q } (u,h),)
\end{align}
where  $u_*(x,t)= \sin (t)$
and consider the equation
\begin{align}\label{eq to bifurcate}
F_\lambda(u,h)=0.
\end{align}
By construction, if $F_\lambda(u,h)=0$, then the function $\tilde u = u_* + u$ solves the problem \eqref{eq:Proe11}. 

Moreover, we have
$$
F_\lambda(0,0) = 0 \qquad \text{for all $\lambda >0$.}
$$
To check if this trivial branch of solutions admits bifurcation, we need to consider the derivative $D F_\lambda(0,0)$.
By direct computation, the operator  $D F_\lambda(0,0)$  is given by
\begin{align}\label{eq diffope}
D F_\lambda(0,0)(v,g)=  \Bigl(L_\lambda^1 v + \cL_\lambda^g u_*, v_\tau (\cdot, (2 m+1) \pi)-g(\cdot)  \Bigl),
\end{align}
with 
$$
\cL_\lambda^g:= \Bigl(\partial_h \Big|_{h=1}L_\lambda^h\Bigr)g =  2g(x)\partial_{tt}+2\lambda g'(x)t \partial_{tx}+ \lambda g''(x)t \partial_t.
$$
For a function $w$ depending only on $t$, we have
$$
\cL_\lambda^g w  = 2g(x)\partial_{tt} + \lambda g''(x)t w_t.
$$
Consequently, the first coordinate in $D F_\lambda(0,0)(v,g)$ is given by 
$$
L_\lambda^1 v + \cL_\lambda^g u_* = v+ \lambda v_{xx} + v_{tt} -2g(x) \sin t + \lambda g''(x)t \cos t,
$$
while the second is 
$$v_\tau (\cdot, (2 m+1) \pi)-g(\cdot).$$
We set 
\begin{align}\label{eqinroU}
U(x,t):=v(x,t)+ g(x)t \cos (t).
\end{align}
Then
we have by direct computation
\begin{align}\label{eequivla}
L_\lambda^1 v + \cL_\lambda^g u_* &= v+ \lambda v_{xx} + v_{tt} -2g(x) \sin t + \lambda g''(x)t \cos t\nonumber\\
&= U+ \lambda U_{xx} + U_{tt}\nonumber\\
v_\tau (\cdot, (2 m+1) \pi)-g(\cdot)&=U_\tau (\cdot, (2 m+1) \pi)
\end{align}
  and 
\begin{align}\label{eq diffope2}
D F_\lambda(0,0)(v,g)=  \Bigl( U+ \lambda U_{xx} + U_{tt}, U_\tau (\cdot, (2 m+1) \pi)  \Bigl).
\end{align}

\section{Computation of the kernel of $D F_\lambda(0,0)$}\label{exkus}
In this section, we  determine the kernel of the operator $D F_\lambda(0,0)$. \\

So suppose that $(v,g)$ is an element of the kernel. Then from \eqref{eq diffope} and  \eqref{eequivla}, the function $U$  is an odd function in $t$ and $2\pi$ periodic  and even in $x$ satisfying
\begin{align}\label{eq. for U}
 U+ \lambda U_{xx} + U_{tt}=0
\end{align}
together with the conditions
\begin{align}\label{eqbounU}
U_\tau (\cdot, (2 m+1) \pi)=0,\quad U (\cdot, \pm(2 m+1) \pi))= \mp (2 m+1) \pi)g(\cdot).
\end{align}

In what follows, we  consider the Fourier coefficients
\begin{align}\label{eqfouriercoefs}
U_k(t):= \frac{1}{\sqrt{2\pi}}\int_{0}^{2\pi} U(x,t)\cos (kx)\,dx, \qquad
g_k:= \frac{1}{\sqrt{2\pi}} \int_{0}^{2\pi} g(x)\cos (kx)\,dx.
\end{align}

Multiplying  \eqref{eq. for U}  and \eqref{eqbounU} with $\cos (kx)$ and integrating in the $x$-variable from $0$ to $2\pi$, $U_k$ is an odd function of class $C^{2,\alpha}([-(2m+1)\pi,(2m+1)\pi])$ solving the initial value problem
\begin{align}\label{initialfor U}
U_k'' + (1-k^2 \lambda) U_k =0
\end{align}
and 
\begin{align}\label{boudaries for U}
U'_k(\pm (2m+1)\pi)&=0 \quad \text{and} \quad U_k(\pm (2m+1)\pi)=\mp (2 m+1) \pi)g_k.
\end{align}

\subsection{Exkursion on the  initial value problem \eqref{initialfor U}-\eqref{boudaries for U}}

\subsection*{In the case  $g_k=0$}
If $g_k=0$, then  $U_k$ solves 
$$
U_k'' + (1-k^2 \lambda) U_k  = 0, \qquad U_k((2 m+1)\pi)=0=U'_k(\pm (2m+1)\pi) 
$$
and we see that $U_k=0$. In this case, we see with \eqref{eqinroU} that the kernel of $D F_\lambda(0,0)$ is trivial. 

\subsection*{In the case  $g_k \not = 0$}

If $g_k \not = 0$, we may assume $g_k = 1$ by considering $U_k /g_k$ in place of $U_k$, so we check whether the linear initial value problem 
\begin{align}\label{eq:v-k-initial-value}
&U_k'' + (1-k^2 \lambda) U_k  = 0,\nonumber\\
 &U'_k(\pm (2m+1)\pi)=0 \quad \text{and} \quad   U_k(\pm (2m+1)\pi)=\mp (2 m+1) \pi).
\end{align}
We put $\mu= k^2 \lambda$, and  let $u_\mu: \R \to \R$ denote the unique solution of the initial value problem
\begin{equation}
  \label{eq:v-k-initial-vu}
u_\mu'' + (1-\mu) u_\mu  =0\qquad u_\mu( \pm (2m+1)\pi)= \mp (2m+1)\pi
\end{equation}
together with
$$u'_\mu(\pm (2m+1)\pi)=0.$$
Then  
\begin{align}\label{eq:v-k-initial-vurelation}
u_\mu(t)\equiv  u (\mu, t)\quad \textrm{and}\quad  U_k(t)= u (k^2 \lambda, t).
\end{align}

\subsection*{The case $\mu =1$}

In this case the \eqref{eq:v-k-initial-vu} reads 
$$
u_\mu'' =0\qquad u_\mu( \pm (2m+1)\pi)= \mp (2m+1)\pi,
$$
together with
$$u'_\mu(\pm (2m+1)\pi)=0.$$
A fundamental system of the linear equation is given by
$$
\phi_1^\mu(t)= 1 , \qquad  \phi_2^\mu(t)=t,
$$
and we have $u_\mu(t) =-t$ and  $u'_\mu(\mp (2m+1)\pi)) =-1< 0$. Therefore, $u_\mu$ is not solution to  \eqref{eq:v-k-initial-vu}.

\subsection*{The case $0 \le \mu < 1$.}

In this case, a fundamental system of the linear equation is given by
$$
\phi_1^\mu(t)= \cos (\sqrt{1-\mu}t), \qquad  \phi_2^\mu(t)= \sin (\sqrt{1-\mu}t)
$$
and 
$$
u_\mu (t) =A\cos (\sqrt{1-\mu}t)+B \sin (\sqrt{1-\mu}t)
$$
for some  real constants $A$ and $B$.
Since we are looking  for an odd function, we must have $u_\mu (0)=0$, which implies $A=0$.
Hence 
\begin{equation}
\label{eq:v-k-change solu22}
u_\mu (t) =B \sin (\sqrt{1-\mu}t).
\end{equation}

In the special case $\mu = 0$, we get

$$
u_\mu(t) = B \sin t.$$
But then 
$
u_\mu ((\mp (2m+1)\pi))) =0\ne \mp (2m+1)\pi,
$
meaning we have no solution in this case.\\

In the general case $\mu \in (0,1)$,\\

If  $\sqrt{1-\mu}(2m+1)$  is an integer,  we see with \eqref{eq:v-k-change solu22} that 

$
u_\mu ((\mp (2m+1)\pi))) =0\ne \mp (2m+1)\pi
$
and we  have no solution.\\
 
We therefore assume that $\sqrt{1-\mu}(2m+1)$  is  not an integer. Then the conditions   $u_\mu( \pm (2m+1)\pi)= \mp (2m+1)\pi$ with \eqref{eq:v-k-change solu22} yield 
$$
u_\mu(t)= -\frac{(2m+1)\pi}{\sin (\sqrt{1-\mu}(2m+1)\pi)} \sin (\sqrt{1-\mu}t)
$$
and 
$$
u'_\mu((2m+1)\pi)= -\frac{(2m+1)\pi\sqrt{1-\mu}}{\sin (\sqrt{1-\mu}(2m+1)\pi)} \cos (\sqrt{1-\mu}(2m+1)\pi).
$$
Consequently 
$$
u'_\mu((2m+1)\pi)=0 
$$
iff and only if  
\begin{align}\label{seqofzeroes}
\sqrt{1-\mu}(2m+1)= \frac{1}{2}+\ell, \quad \ell \in \N.
\end{align}
Since   $\mu \in (0,1)$, $ \sqrt{1-\mu}< 1$ and hence $ 0\leq \ell \leq 2m$.

For $m=0$ and  $\mu \in (0,1)$, $\sqrt{1-\mu}$ is not an integer and $\sqrt{1-\mu}=\frac{1}{2}+\ell$ if and only if 
$\ell=0$ and  
\begin{align}\label{eqmequqls0}
\mu=\frac{3}{4}.
\end{align}

\subsection{The case $\mu>1$.}

A fundamental system of the linear equation is then given by
$$
\phi_1^\mu(t)= e^{-\sqrt{\mu -1}t}, \qquad  \phi_2^\mu(t)= e^{\sqrt{\mu -1}t}
$$
and 
 \begin{equation}
\label{eq:v-k-change solution3}
u_\mu(t)=-\frac{(2m+1)\pi}{\sinh  (\sqrt{1-\mu}(2m+1)\pi)} \sinh (\sqrt{1-\mu}t) 
\end{equation}
and
 \begin{equation}
\label{eq:v-k-change solution3}
u'_\mu(t)=-\frac{(2m+1)\pi\sqrt{1-\mu}}{\sinh  (\sqrt{1-\mu}(2m+1)\pi)} \cosh (\sqrt{1-\mu} (2m+1)\pi)) <0.\\\\
\end{equation}

\subsection{Computation of the kernel}\label{seckernel}
Summarizing what we have got so far, for every $m\in \mathbb{N}$, \eqref{seqofzeroes} yields strictly  decreasing sequence of $(\mu_\ell(m))_{0\leq \ell \leq 2m}$ made of zeroes of the function 
\begin{align}\label{eq:grapheigen}
V: \mu \longmapsto u'_\mu((2m+1)\pi),
\end{align}
and such that 
\begin{align}\label{eq:foroodin1}
u_{ \ell}(t)= -(-1)^\ell (2m+1)\pi  \sin (\sqrt{1-\mu_\ell(m)}t)
\end{align}
is solution to  \eqref{eq:v-k-initial-vu}. 

Moreover, for each  ${0\leq \ell \leq 2m}$,  $ \mu_\ell(m)\in (0, 1)$ and 
\begin{align}\label{seqofzeroes2}
 \mu_\ell(m)&=1-\frac{1}{4}\Bigl(   \frac{1+2\ell }{1+2m}\Bigl)^2=\frac{\Bigl(1+2(2m-\ell))\bigl)\Bigl(3+2(2m+\ell))\bigl)}{4\Bigl(1+2m \Bigl)^2}.
\end{align}

Recalling   \eqref{eq:v-k-initial-vurelation}  it follows that for  a suitable  $\lambda>0$, the kernel of the operator $D F_\lambda(0,0)$ is determined by finding the functions $ U_k$, where the  integers $k$  satisfy the restriction
\begin{align}\label{eq:foroofc1}
0<k^2 \lambda<1.
\end{align}
Together with 
\begin{align}\label{eq:foroofc21}
k^2 \lambda =  \mu_\ell(m), \quad \textrm{for some }\quad {0\leq \ell \leq 2m}
\end{align}
Now  choose  any  $0 \leq \ell_0 \leq 2m$ and take  $\lambda=  \mu_{\ell_0}(m)$. Then  \eqref{eq:foroofc21} we are led to finding the $k$'s such that
\begin{align}\label{eq:foroofc22}
k^2 =  \frac{\mu_\ell(m)}{\mu_{\ell_0}(m)}=\frac{\Bigl(1+2(2m-\ell)\bigl)\Bigl(3+2(2m+\ell)\bigl)}{\Bigl(1+2(2m-\ell_0)\bigl)\Bigl(3+2(2m+\ell_0))\bigl)}
 \quad \textrm{for some }\quad {0\leq \ell \leq 2m}.
\end{align}
Since the sequence  $(\mu_\ell(m))_{0\leq \ell \leq 2m}$  is strictly  decreasing, 
$\frac{\mu_\ell(m)}{\mu_{\ell_0}(m)}< 1 \quad \textrm{for every }\quad {\ell >\ell_0}$.
So we might assume in  \eqref{eq:foroofc22} that  $\ell  \leqslant \ell_0$.\\
 
Obviously $k=1$  for  $\ell  =\ell_0$.

For  $\ell  < \ell_0$, we check that the ratio in the right hand side of  \eqref{eq:foroofc22} is not an integer. Indeed, if there exists some   $\ell  < \ell_0$ such that $$\frac{\Bigl(1+2(2m-\ell)\bigl)\Bigl(3+2(2m+\ell)\Bigl)}{\Bigl(1+2(2m-\ell_0)\bigl)\Bigl(3+2(2m+\ell_0)\Bigl)}\in \mathbb{N}.$$
Then $(3+2(2m+\ell_0))$ must divide at least  one of the factors  $(1+2(2m-\ell))$  or $(1+2(2m-\ell_0))$. \\

 If  $(3+2(2m+\ell_0))$ divides  $(1+2(2m-\ell))$, necessarily   $(3+2(2m+\ell_0))\leqslant  (1+2(2m-\ell))$. That is $2+2\ell_0\leqslant  -2\ell$.
 
If  $(3+2(2m+\ell_0))$ divides  $(3+2(2m+\ell))$, the same argument  implies the inequality  $ \ell_0 \leqslant \ell$, which is in contradiction with that assumption $\ell  < \ell_0$.\\
 
We  can therefore derive from \eqref{eqinroU} and  \eqref{eq:foroodin1} the that  for all $0 \leq \ell \leq  2m$, the kernel of  the operator  $D F_{\mu_{\ell}}(0,0)$ is one-dimensional and spanned by $(v_\ell,g_\ell)$ with
$$
v_\ell (x,t)=  w_{\mu_\ell(m)}(t)\cos x, \qquad g_\ell(x)=\cos x,
$$
with
$$w_{\mu_\ell(m)}(t)= -(-1)^\ell(2m+1)\pi \sin (\sqrt{1-\mu_\ell(m)}t)-t\cos (t)
$$
and 
$$\sqrt{1-\mu_\ell(m)}(2m+1)=\frac{1}{2}+\ell.$$

We underline  that the particular case  $m=0$ gives    $\mu_0=\frac{3}{4}$ and  the  kernel  of $D F_{\mu_0}(0,0)$  spanned by $(v_0,g_0)$ with
$$
v_0(x,t)=  w_{\mu_*}(t)\cos x, \qquad g_0(x)=\cos x,
$$   
with
$$
w_{\mu_*}(t)= -\pi \sin (t/2)-t\cos (t).
$$

We  also observe  from \eqref{seqofzeroes2} that   $\mu_m(m)=\frac{3}{4}$ and therefore,
\begin{align}\label{eq:forfixed m}
\mu_\ell(m)  \geqslant \mu_{m}(m) =\frac{3}{4} \quad \textrm{ for all}\quad  0 \leq \ell \leq m.
\end{align}

In the next section, we  gather the required assumptions  to apply the  bifurcation result from simple eigenvalues of  Crandall-Rabinowitz to prove existence of branches of solution to the equation \eqref{eq to bifurcate}.

\section{Properties of the linearized operator}\label{Spectral}
In this section,  we analyse the operator $DF_{\mu_\lambda}(0,0)$ and gather the assumptions that enable us to apply  the  Crandall-Rabinowitz bifurcation theorem.  \\

In the following, we  choose $\ell$ between  $0$ and $2m$  and  let $$\ker(DF_{\mu_\ell(m)}(0,0))^{\perp}\subseteq X_2\times Y_2\subseteq L^2_{p,e}((-\pi,\pi)\times  (-(2 m+1) \pi, (2 m+1) \pi) )$$ be the complement of $DF_{\mu_\ell(m)}(0,0)$ with respect to the scalar product  
\nolinebreak
\begin{equation}\label{scalar-product}
\langle(v,g),(w,h)\rangle:=\int_{(-\pi,\pi)\times  (-(2 m+1) \pi, (2 m+1) \pi) }v(x,t)w(x,t)\,dxdt+\int_{(-\pi,\pi)}g(x)h(x)dx.
\end{equation} 
Next we consider the restriction mapping
$$
DF_{\mu_\ell(m)}(0,0): \ker(DF_{\mu_\ell(m)}(0,0))^{\perp}\ni(w,h)\mapsto DF_{\mu_\ell(m)}(0,0) \cdot\left[ {\begin{array}{cc}
   w \\
   h\\
  \end{array} } \right] \in X_0 \times Z_2
$$
and set 
\begin{align}\label{MatriF}
 A(w)(x):&= \int_{(-\pi_m,\pi_m)}t\cos(t)\partial_{xx}w(x,t)\,dt, \quad &B(w)(x): = \int_{(-\pi_m,\pi_m)}\sin(t)w(x,t)\,dt \nonumber\\
  \mathcal{L}(z)\cdot v:&= \int_{(-\pi,\pi)}v_t (x,\pi_m) z(x)\,dx,  \quad & \mathcal{C}(w) \cdot v:= 2\int_{(-\pi,\pi)}v_t (x,\pi_m) w(x,\pi_m)\,dx \nonumber\\
 \mathcal{K}(z)\cdot g:& = -\int_{(-\pi,\pi)}g(x)z(x)\,dx. 
  \end{align}
With these notations, we have
\begin{Lemma}\label{lemorthoKer}
Let $(v,g)\in \ker(DF_{\mu_\ell(m)}(0,0))^{\perp}$ and $(w,z)\in  X_0\times Z_2$. Then,
\begin{align}\label{op:adjoint}
   \Big\langle DF_{\lambda}(0,0)\left[ {\begin{array}{cc}
   v \\
   g\\
  \end{array} } \right],\left[ {\begin{array}{ccc}
   w \\
   z\\
  \end{array} } \right]\Big\rangle 
  &=\Big\langle\left[ {\begin{array}{ccc}
   \id+\lambda \partial_{xx}+\partial_{tt}+ \mathcal{C}  & \calL  \\
   \lambda A-2B &   \mathcal{K} \\
  \end{array} } \right]\left[ {\begin{array}{cc}
   w \\
   z\\
  \end{array} } \right],\left[ {\begin{array}{cc}
   v \\
   g\\
  \end{array} } \right]\Big\rangle.
  \end{align}
  
Moreover for every  \textrm{$ 0\leq \ell \leq 2m$}, we have
\begin{align}\label{decompol1}
X_0\times Z_2= DF_{\mu_\ell(m)}(0,0)\Big(\ker(DF_{\mu_\ell(m)}(0,0))^{\perp}\Big)\oplus E_\ell, 
\end{align}
where  
\be 
E_\ell := \spann\big(\ov w^\ell,\ov z^\ell \big)=\{a(\ov w^\ell,\ov z^\ell): a\in\R\},\ee
and 
\begin{align}\label{eqdisolcokernell}
\ov w^\ell(x,t):= \sin (\sqrt{1-\mu_\ell(m)}t)\cos (x)\quad \textrm{and} \quad \ov z^\ell(x):= -2(-1)^\ell\cos (x).
\end{align}
\end{Lemma}

\proof
Let $(v,g)\in \ker(DF_{\mu_\ell(m)}(0,0))^{\perp}$ and $(w,z)\in  X_0\times Z_2$. To shorten we let $\pi_m:= (2 m+1) \pi$.  Then 
  \begin{align*}
 & \Big\langle DF_{\lambda}(0,0)\left[ {\begin{array}{cc}
   v \\
   g\\
  \end{array} } \right],\left[ {\begin{array}{cc}
   w \\
   z\\ 
  \end{array} } \right]\Big\rangle\nonumber\\
  &=\int_{(-\pi,\pi)\times (-\pi_m,\pi_m) }\big(v+\lambda v_{xx}+v_{tt}-2g(x)\sin(t)+\lambda g''(x)t\cos(t)\big)w(x,t)\,dxdt\nonumber\\
  &+\int_{(-\pi,\pi)}\bigl(v_t(x, \pi_m)-g(x)\bigl)z(x)\,dx
  \end{align*}
  Performing an integration by parts, we find 
  \begin{align*}
   &\Big\langle DF_{\lambda}(0,0)\left[ {\begin{array}{cc}
   v \\
   g\\
  \end{array} } \right],\left[ {\begin{array}{ccc}
   w \\
   z\\
  \end{array} } \right]\Big\rangle \\
  &= \int_{(-\pi,\pi)\times (-\pi_m,\pi_m) }w(v+\lambda  v_{xx}+v_{tt})\,dxdt+\l\int_{(-\pi,\pi)\times (-\pi_m,\pi_m) }t\cos(t)g(x)w_{xx}(x,t)\,dxdt\\
  &-2\int_{(-\pi,\pi)\times (-\pi_m,\pi_m) }\sin(t)g(x)w(x,t)\,dxdt+\lambda \int_{-\pi_m}^{\pi_m} (wv_x-vw_x)\Big|_{-\pi}^{+\pi}dt+\int_{-\pi}^\pi(v_tw-w_tv)\Big|_{-\pi_m}^{+\pi_m}dx\\
  &+\lambda \int_{-\pi_m}^{+\pi_m}t\cos(t)(g'w-gw_x)\Big|_{-\pi}^{+\pi}dt+\int_{(-\pi,\pi)}\bigl(v_t(x, \pi)-g(x)\bigl)z(x)\,dx\\
  &=\int_{(-\pi,\pi)\times (-\pi_m,\pi_m) }v(w+\lambda w_{xx}+w_{tt})\,dxdt+\int_{-\pi}^{\pi}g(x)\big(\l A(w)-2B(w)\big) dx\\
  &+\int_{(-\pi,\pi)}\bigl(v_t(x, \pi_m)-g(x)\bigl)z(x)\,dx+2\int_{-\pi}^\pi v_t (x,\pi_m) w(x,\pi_m) dx+S(w,v,g)\\
  &=\Big\langle\left[ {\begin{array}{ccc}
   \id+\lambda \partial_{xx}+\partial_{tt}+ \mathcal{C}  & \calL  \\
   \lambda A-2B &   \mathcal{K} \\
  \end{array} } \right]\left[ {\begin{array}{cc}
   w \\
   z\\
  \end{array} } \right],\left[ {\begin{array}{cc}
   v \\
   g\\
  \end{array} } \right]\Big\rangle+S(w,v,g)\\
  &=\Big\langle\mathcal{F}_\lambda \left[ {\begin{array}{cc}
   w \\
   z\\
  \end{array} } \right],\left[ {\begin{array}{cc}
   v \\
   g\\
  \end{array} } \right]\Big\rangle+S(w,v,g),
  \end{align*}
  where we have used \eqref{MatriF} and $S(w,v,g)$ is given by
  \begin{align*}
  S(w,v,g):& = \l \int_{-\pi_m}^{\pi_m}(w\partial_x v-v\partial_xw)\Big|_{-\pi}^{+\pi}dt+\l  \int_{-\pi_m}^{\pi_m}t\cos(t)(g'w-gw_x)\Big|_{-\pi}^{+\pi}dt\\
  &-2\int_{(-\pi,\pi)}v(x, \pi_m)w_t (x,\pi_m) \,dx.
     \end{align*}
Since $v=0$ on $\partial\O$, $\partial_xw,\partial_xv$ and $g'$ are $2\pi$ periodic, and $w$ is even in $x$ and $t$,  it follows that $S(w,v,g)=0$. 
Therefore we have 
  \be\label{adjoint-op}
  \Big\langle DF_{\lambda}(0,0)\left[ {\begin{array}{cc}
   v \\
   g\\
  \end{array} } \right],\left[ {\begin{array}{cc}
   w \\
   z\\
  \end{array} } \right]\Big\rangle = \Big\langle \mathcal{F}_\lambda \left[ {\begin{array}{cc}
   w \\
   z\\
  \end{array} } \right],\left[ {\begin{array}{cc}
   v \\
   g\\
  \end{array} } \right]\Big\rangle.
 \ee  
This proves \eqref{op:adjoint}.\\
 
To prove the splitting  \eqref{decompol1}, we let $(\ov w,\ov z)$ such that $ \mathcal{F}_\lambda \left[ {\begin{array}{cc}
   \ov w \\
   \ov z\\
  \end{array} } \right] = 0$, that is 
  \be\label{kernel-adjoint-op}
  \left\{\begin{aligned}
  \ov w+\l \partial_{xx}\ov w+\partial_{tt}\ov w+ \mathcal{C}(\ov w)+ \calL(\ov z&)=0 \\
\l A(\ov w)-2B(\ov w)+\mathcal{K}(\ov z&)= 0.
 \end{aligned}
\right.
  \ee
We claim that 
\be\label{Claim B8}
\textit{The PDEs in  \eqref{kernel-adjoint-op} is uniquely solvable in $X_2\times Y_2$ (up to a multiplicative constant).}
\ee
Indeed writing  $\ov w$ and $\ov z$ in the form 
  $$
  \ov w = \sum_{k\geq 0}\ov w_k(t)\cos(kx)\quad\text{and}\quad \ov z = \sum_{k\geq 0}\ov z_k\cos(kx)
  $$
the second equation in \eqref{kernel-adjoint-op} reads 
 \be\label{B18}
 \l \int_{-\pi_m}^{\pi_m}t\cos(t)\partial_{xx}\ov w(x,t)\,dt-2\int_{-\pi_m}^{\pi_m}\sin(t)\ov w(x,t)\,dt+\mathcal{K}(\ov z) = 0,\qquad x\in(-\pi,\pi).
 \ee 
Next we apply  $g(x)=\cos(kx)$  to \eqref{B18} and integrate  (by using Fubini Theorem) to get
 \be\label{B19}
 \int_{-\pi_m}^{\pi_m}\big(k^2\l t\cos(t)+2\sin(t)\big)\ov w_k(t)\,dt+\ov z_k=0.
 \ee
Similarly  we apply the first equation in \eqref{kernel-adjoint-op}  to $v(x,t)=\cos(kx)$  and integrate  to  find 
 \be\label{B20}
 \ov w_k''(t)+(1-k^2\l )\ov w_k(t)=0.
 \ee
Next we apply  it again  to  $v(x,t)=\sin(t)\cos(kx)$ and derive
\be \label{Bz0}
\ov z_k =-2\ov w_k (\pi_m).
\ee
We set $\mu:= k^2\l=  k^2\mu_\ell(m)$. Then 

If $\mu =1$,  then the solution $\ov w(t)$  to the equation 
\be\label{B200}
 \ov w ''(t)+(1-\mu )\ov w(t)=0.
 \ee
is given by $\ov w=A t$ and integrating by parts, we see with \eqref{Bz0} that $A=0$. \\

If  $0 \leqslant \mu < 1$, the general solution of the  equation  \eqref{B20} is given by
\begin{align}\label{C444}
\ov w (t)=B \sin (\sqrt{1-\mu}t),
\end{align}
for some real constant $B$ and we deduce 
\begin{align}\label{C45}
\ov z_k =-2 B \sin (\sqrt{1-\mu }\pi_m). 
\end{align}
To  check if $\ov w$ and $ \ov z$ satisfy \eqref{B19}, we rewrite this  condition in term of $\mu$.\\ 

We set $\xi:=\sqrt{1-\mu}$.
Then a straightforward computation allows to get 
\begin{align}\label{C23}
\int_{0}^{\pi_m} \sin(t) \sin(\xi t)  \,dt&= \frac{1}{2}\int_{0}^{\pi_m} \Bigl( \cos((1-\xi)t) - \cos((1+\xi)t) \,dt=\frac{1}{\mu} \sin(\xi \pi_m).
\end{align}
We also have 
\begin{align}\label{C22}
\int_{0}^{\pi_m}\big(\mu t\cos(t) \sin(\xi t) \,dt&=\frac{\mu}{2}\int_{0}^{\pi_m} t\Bigl( \sin((1+\xi)t) - \sin((1-\xi)t)\Bigl) \,dt\nonumber\\
&=\frac{\mu \pi_m}{2} U(\pi_m)-\frac{2-\mu}{\mu} \sin(\xi \pi_m),
\end{align}
where the last line follows after integrating by parts and 
$$ U(t)=-\frac{1}{1+\xi} \cos((1+\xi)t) + \frac{1}{1-\xi} \cos((1-\xi)t)\quad \textrm{and}\quad U(\pi_m)= -\frac{2\sqrt{1-\mu} }{\mu}\cos(\sqrt{1-\mu} \pi_m).$$
Gathering \eqref{C23} and  \eqref{C22}, it follows that
\be\label{B09}
  \int_{-\pi_m}^{\pi_m}\big(\mu t\cos(t)+2\sin(t)\big)\sin(\xi t)\,dt = 2\sin(\sqrt{1-\mu} \pi_m)-2\sqrt{1-\mu} \pi \cos(\sqrt{1-\mu} \pi_m).
 \ee
This with \eqref{C45} give
\be\label{Bverified}
  B \int_{-\pi_m}^{\pi_m} \big(\mu t\cos(t)+2\sin(t)\big)\sin(\xi t)\,dt+\ov z_k  = -2B\sqrt{1-\mu} \pi_m \cos(\sqrt{1-\mu} \pi_m).
 \ee
It is plain  from \eqref{Bverified} that $B=0$ for $\mu=0$. In this case the solution of \eqref{kernel-adjoint-op}  given  \eqref{C444} and \eqref{C45} is trivial. In particular  $\ov z_0=0$ and $\ov w_0=0$. 
Furthermore,  we  see with \eqref{Bverified} that the condition \eqref{B19} is fulfilled  with $B\ne 0$  if and only if 
\be\label{Bverified2}
k^2\mu_\ell(m)=\mu_{\bar{\ell}}(m) \quad\textrm{ for some $0\leqslant \bar{\ell}\leqslant 2m$},
 \ee
which from the analysis in Subsection \ref{seckernel} is equivalent to  $k=1$.\\

To complete the proof, we show that  $\ov z_k=0$ and $\ov w_k=0$ for every  $k\geqslant2$. To see this, we distinguish two cases:\\

If  $k\geqslant2$ is such that  $\mu=k^2\mu_\ell(m)<1$,  then from  the argument following  \eqref{Bverified2},  we must have $B=0$.\\

Now if   $k\geqslant2$  is such that  $\mu=k^2\mu_\ell(m)>1$, then the general solution  of \eqref{B20} in this case is   given by 
$$ \ov w_k(t)=A \bigl( e^{\mu (k) t}-e^{-\mu (k)t}\bigl)=2 A \sinh (\mu (k) t),$$
where 
$$\mu (k):=\sqrt{k^2 \mu_{\ell}(m)-1}.$$
Using \cite[Page 231, $1.$]{I.S.I.M}, we find
\begin{align}\label{B22}
&\int_{0}^{\pi_m}\sin(t) \sinh (\mu (k) t)\,dt =\frac{1}{1+\mu (k)^2} \sinh (\mu (k)  \pi_m).
\end{align} 
Also \cite[Page 230, $6.$]{I.S.I.M} yields,
\begin{align}
\int_{0}^{\pi_m}\big(t\cos(t)  e^{\mu (k) t}\,dt &=\frac{1}{1+\mu (k)^2}\Bigl[e^{\mu (k)\pi_m} \Bigl( \frac{\mu (k)^2-1}{1+\mu (k)^2} -\mu (k) \pi_m \Bigl)+\frac{\mu (k)^2-1}{1+\mu (k)^2}\Bigl],\nonumber\\
\int_{0}^{\pi_m}\big(t\cos(t)  e^{-\mu (k) t}\,dt &=\frac{1}{1+\mu (k)^2}\Bigl[e^{-\mu (k)\pi_m} \Bigl( \frac{\mu (k)^2-1}{1+\mu (k)^2} +\mu (k)  \pi_m \Bigl)+\frac{\mu (k)^2-1}{1+\mu (k)^2}\Bigl],\nonumber
\end{align} 
so that 
\begin{align}
\int_{0}^{\pi_m} \mu t \cos(t) \sinh (\mu (k) t)\,dt =\frac{\mu}{1+\mu (k)^2}\Bigl[ \Bigl( \frac{\mu (k)^2-1}{1+\mu (k)^2} \Bigl) \sinh (\mu (k) \pi_m) -\mu (k) \pi_m  \cosh (\mu (k) \pi_m)\Bigl)\Bigl],
\end{align} 
and  
\begin{align}\label{B32}
&\int_{-\pi_m}^{\pi_m}\big(\mu t\cos(t)+2\sin(t)\big) \sinh (\mu (k) t)\,dt\nonumber\\
&=\frac{2}{1+\mu (k)^2}\Bigl[ \mu\Bigl( \frac{\mu (k)^2-1}{1+\mu (k)^2} \Bigl)+ 2 \Bigl) \sinh (\mu (k) \pi_m) -\mu \mu (k) \pi_m  \cosh (\mu (k) \pi_m)\Bigl)\Bigl]\nonumber\\
&=2\Bigl(\sinh (\mu (k) \pi_m)-\mu (k) \pi_m  \cosh (\mu (k) \pi_m)\Bigl).
\end{align} 
From this
\begin{align}\label{B32}
\int_{-\pi_m}^{\pi_m}\big(\mu t\cos(t)+2\sin(t)\big) \ov w_k(t)&=2A \int_{-\pi_m}^{\pi_m}\big(\mu t\cos(t)+2\sin(t)\big)\sinh (\mu (k) t)\,dt+\ov z_k \nonumber\\
&=4A \Bigl(\sinh (\mu (k) \pi_m)-\mu (k) \pi_m  \cosh (\mu (k) \pi_m )\Bigl)-4A \sinh (\mu (k) \pi_m)\nonumber\\
&=-4A\mu (k) \pi_m  \cosh (\mu (k) \pi_m)
\end{align} 
and we deduce that $A=0.$ 

Next, if  $k\geqslant2$ is such that  $\mu=k^2\mu_\ell(m)=1$, we already know from the argument following \eqref{B200} that  the corresponding $\ov w_k$ vanishes. Gathering  all the steps, we get  $\ov z_k=0$ and $\ov w_k=0$  for every  $k\geqslant2$.\\

Finally, the solution $\ov w$ to  \eqref{kernel-adjoint-op} is given by \begin{align}\label{eqdisolcokernell}
\ov w^\ell(x,t):= \sin (\sqrt{1-\mu_\ell(m)}t)\cos (x)\quad \textrm{and} \quad \ov z^\ell(x):= -2(-1)^\ell\cos (x)
\end{align}
and setting 
\be 
E_\ell := \spann\big(\ov w^\ell,\ov z^\ell \big)=\{a(\ov w^\ell,\ov z^\ell): a\in\R\},\ee
we obtain \eqref{decompol1}. \QED

\begin{Proposition}\label{propPhlinearop}
Let  \textrm{$ 0\leq \ell \leq 2m$}. Then the  linear
operator
$$ 
\cH_\ell:=DF_{\mu_\ell(m)}(0,0): X_2 \times Y_2 \to X_0  \times Z_2
$$
has the following properties:
\begin{itemize}
\item[(i)] The kernel $N(\cH_\ell)$ of $\cH_\ell$ is spanned $(v_\ell,g_\ell)$ with
\begin{equation}\label{eq:kernel}
v_\ell (x,t)=  w_{\mu_\ell(m)}(t)\cos x, \qquad g_\ell(x)=\cos x,
\end{equation}
with
$$w_{\mu_\ell(m)}(t)= -(-1)^\ell(2m+1)\pi \sin (\sqrt{1-\mu_\ell(m)}t)-t\cos (t)
$$
and 
$$\sqrt{1-\mu_\ell(m)}(2m+1)=\frac{1}{2}+\ell.$$
\item[(ii)] The range of $\cH_\ell$ is given by
$$
R(\cH_\ell)= \left\lbrace (v,g)\in X_0\times Z_2: \int_{(-\pi,\pi)\times(-\pi_m,\pi_m) }\ov w_\ell(x,t)v(x,t)\,dxdt+\int_{(-\pi,\pi)}\ov z_\ell(x)g(x)dx=0\right\rbrace,
$$
\end{itemize}
with 
\begin{equation}\label{eq:for image}
\ov w_\ell(x,t):=\sin (\sqrt{1-\mu_\ell(m)}t)\cos (x)\quad \textrm{and} \quad \ov z_\ell(x):=-2(-1)^\ell\cos (x).
\end{equation}
Moreover,
\begin{equation}
  \label{eq:transversality-cond4}
\partial_\lambda \Bigl|_{\lambda=\mu_\ell(m) }DF_{\l}(0,0)(v_\ell,g_\ell)\not  \in \; R(\cH_\ell).\\
\end{equation}
\end{Proposition}

The remaining part of the paper is to prove Proposition \ref{propPhlinearop}. \\

\proof

The proof of $(i)$ is already done at the end of  Section  \ref{exkus}.  Moreover, the mapping $DF_{\mu_\ell(m) }(0,0)$ is  injective. 

To prove $(ii)$ we first observe  from \eqref{decompol1} in Lemma \ref{lemorthoKer} that 
\be
Im\Big(DF_{\mu_\ell(m) }(0,0)\Big)\subseteq E_\ell^{\perp},
\ee
where 
\be
E_\ell^{\perp} := \left\lbrace (v,g)\in X_0\times Y_2: \int_{(-\pi,\pi)\times (-\pi_m,\pi_m) }\ov w^\ell(x,t)v(x,t)\,dxdt+\int_{(-\pi,\pi)}\ov z^\ell(x)g(x)dx=0\right\rbrace.
\ee

Now we will prove that the mapping 
\be
DF_{\mu_\ell(m)}(0,0): \ker(DF_{\mu_\ell(m)}(0,0))^{\perp}\longrightarrow E_\ell^{\perp}
\ee is surjective. \\
 
By Fredholm alternative theorem, for any  $(v,z)\in  X_0\times Z_2$, the PDEs
\be\label{eq-solvability}
DF_{\mu_\ell(m)}(0,0)\cdot\left[ {\begin{array}{cc}
   w \\
   h\\
  \end{array} } \right] = (v,z)
\ee
is solvable in $L^2_{p,e}(-\pi,\pi)\times L^2_{p,e}(-\pi_m,\pi_m)$ in the sense of distribution.  Moreover by using Fourier method, it is not also difficult to show that if $v$ and $g$ are both periodic, then one can find a solution $(w,h)$ that is also periodic. By elliptic regularity we have $(w,h)\in C^{2,\al}_{p,e}(\ov\O_*)\times C^{2,\al}_{p,e}(\R)$. Indeed, letting $$U(x,t):=w(x,t)+h(x)t\cos(t),$$  a direct calculation gives 
\be
DF_{\mu_\ell(m)}(0,0)\cdot\left[ {\begin{array}{cc}
   w \\
   h\\
  \end{array} } \right] = (v,z)\Longleftrightarrow \left\{\begin{aligned}\label{eqregul2}
  \calL_{\mu_\ell(m)}U+U&=v \quad\text{in}\quad \O_*\\
U&= \mp \pi_m h \quad\text{on}\quad \partial\O_*^{\mp}\\
\partial_\nu U &= z \quad\text{on}\quad \partial\O_*,
 \end{aligned}
\right.
\ee
with $\calL_{\mu_\ell(m)}=\mu_\ell(m)\partial_{xx}+\partial_{tt}$. We now have $\partial_\nu U = z\in C^{1,\al}(\R)$  and   by classical elliptic regularity, we deduce that $ U\in C^{2,\al}(\ov \O_*)$.
There result  from $U= \mp \pi_m h$ on $\partial\O_*^{\mp}$  that  $ h \in C^{2,\al}(\R)$ and  also  $ w\in C^{2,\al}(\ov \O_*)$. Finaly, we observe  that
the functions $(x, t)\mapsto -U(x,-t)$  and  $(x, t)\mapsto  \frac{1}{2}\left\lbrace U(x,t)+U(-x,t) \right\rbrace$ also solve the equation  \eqref{eqregul2}. In particluar $U$  is odd in $t$, which implies that $w$ is odd in $t$. In addition, $U$ in $x$ and we deduce  from $U= \mp \pi_m h$ on $\partial\O_*^{\mp}$  that $h$ and $w$ are even in $x$. \\

We complete the proof, we have to check \eqref{eq:transversality-cond4}. That is 
\be
\frac{d}{d\l}_{\Big|\l=\mu_\ell(m)}DF_\l(0,0)\left[ {\begin{array}{cc}
   v_\ell \\
   g_\ell\\
  \end{array} } \right]\notin E_\ell^\perp,
\ee
$$
v_\ell (x,t)=  w_{\mu_\ell(m)}(t)\cos x, \qquad g_\ell(x)=\cos x,
$$
with
$$w_{\mu_\ell(m)}(t)= -(-1)^\ell(2m+1)\pi \sin (\sqrt{1-\mu_\ell(m)}t)-t\cos (t)
$$
and 
$$\sqrt{1-\mu_\ell(m)}(2m+1)=\frac{1}{2}+\ell.$$
Following \eqref{eq diffope2}, we  set 
$$
U_\ell(x,t):=v_\ell(x,t)+ g_\ell(x)t \cos (t)= -(-1)^\ell\pi_m \sin (\sqrt{1-\mu_\ell(m)}t)\cos x.
$$
Then
$$
D F_{\lambda}(0,0)(v_\ell, g_\ell)=  \Bigl((\lambda-\mu_\ell(m))(-1)^\ell\pi_m \sin (\sqrt{1-\mu_\ell(m)}t)\cos x, 0 \Bigl)
$$
and 
\begin{align}
&\Big<(\ov w^\ell,\ov z^\ell),\Bigl((-1)^\ell\pi_m \sin (\sqrt{1-\mu_\ell(m)}t)\cos x, 0 \Bigl)\Big>_{L^2((-\pi,\pi)\times (-\pi_m,\pi_m)} \nonumber\\
&= (-1)^\ell\pi_m \int_{(-\pi,\pi)\times (-\pi_m,\pi_m) }\cos^2(x) \sin^2 (\sqrt{1-\mu_\ell(m)}t)\,dxdt= (-1)^\ell \pi \pi^2_m\ne 0.
\end{align}
\QED

\section{Solving problem \eqref{eq:Proe11}  }\label{eq:ProofTheo1}
In this section we complte the proof of   Theorem \ref{Strip1}.

\noindent {\bf Proof of Theorem \ref{Strip1} (completed).} 
Recalling \eqref{eq:maptFF}, the proof of Theorem \ref{Strip1} will be established by applying 
the Crandall-Rabinowitz Bifurcation Theorem to solve the equation
\begin{align}\label{eq:maptFsolvF}
 F_\lambda(u,h)= (L_\lambda^{1+h} (u+u_*), \widetilde{ Q } (u,h))=(0,0).
\end{align}
 where   $(u,h)\in  X_2 \times Y_2^+$ and   $u_*(x,t)= \sin (t)$.
As already observed in Section \ref{eq.setting},  if \eqref{eq:maptFsolvF} holds, then the function $\tilde u = u_* + u$ solves the problem $(Q)$ in  \eqref{eq:Proe11}.\\
To solve equation \eqref{eq:maptFsolvF}, we  fix  $\ell $ betweeen $ 0$ and $2m$, and consider the operator
$$\cH_\ell:=DF_{\mu_\ell(m)}(0,0).$$
 in Proposition \eqref{propPhlinearop},
as well as the space 
\begin{align}\label{eq:spaorthokernel}
\mathcal{X}_\ell^{\perp} := \left\lbrace (v,g)\in X_2 \times Y_2^+: \int_{(-\pi,\pi)\times(-\pi_m,\pi_m) }v(x,t)v_\ell(x,t)\,dxdt+\int_{(-\pi,\pi)} g(x)g_\ell(x)dx=0\right\rbrace.
\end{align}
By Proposition~\eqref{propPhlinearop} and  the Crandall-Rabinowitz Theorem (see \cite[Theorem 1.7]{M.CR}), we then find ${\e_\ell}>0$ and a smooth curve
$$
(-{\e_\ell},{\e_\ell}) \to   \R_+ \times X_2\times  Y_2^+ , \qquad s \mapsto (\lambda^\ell_s,\varphi^\ell_s, \psi^\ell_s)
$$
such that
\begin{enumerate}
\item[(i)] $F_{\lambda_s} (\varphi_s, \psi_s)=0$ for $s \in (-{\e_0},{\e_0})$,
\item[(ii)] $ \lambda^\ell_0= \mu_\ell(m)$, and
\item[(iii)]  $ (\varphi^\ell_s, \psi^\ell_s) = s \bigl((v_\ell,g_\ell) + (\mu^\ell_s, \kappa^\ell_s)\bigr)$ for $s \in (-\e_0,\e_0)$ with a smooth curve
$$
(-{\e_0},{\e_0}) \to \mathcal{X}_\ell^{\perp}, \qquad s \mapsto  (\mu^\ell_s, \kappa^\ell_s)
$$
satisfying $ (\mu^\ell_0, \kappa^\ell_0)=(0,0)$
and
$$\int_{(-\pi,\pi)\times(-\pi_m,\pi_m) } \mu^\ell_s v_\ell(x,t)\,dxdt+\int_{(-\pi,\pi)}  \kappa^\ell_s g_\ell(x)dx=0.$$
\end{enumerate}

Since $(\lambda^\ell_s,\varphi^\ell_s, \psi^\ell_s)$ is a solution to
\eqref{eq:maptFsolvF} for every $s\in (-\e_0,\e_0)$,  the
function $$\tilde u^\ell_s = u_* +  s \bigl(v_\ell + \mu^\ell_s\bigr).$$
 solves the
overdetermined boundary value  problem $(Q)$ in  \eqref{eq:Proe11}.
Recalling the ansatz  \eqref{eqans}, the function
\begin{align}\label{eqans2}
v^\ell_s(x,t)= \tilde u^\ell_s(x ,(1+ s(g_\ell(x) +\kappa^\ell_s(x))t)
\end{align}
solves the problem  (P) in \eqref{eq:perturbed-strip2} on the domain 
$$
\Omega^\ell_s:= \Biggl\{\Bigl(x, \frac{\t}{1+ s(g_\ell(x) +\kappa^\ell_s(x))}\Bigl)\::\: (x,\t) \in  \Omega_*\Biggl\}.
$$
The proof is complete. \QED

\section{Crandall-Rabinowitz bifurcation theorem} \label{eq: Cradal Rabi1}

\begin{Theorem}[Crandall-Rabinowitz bifurcation theorem, \cite{M.CR}] \label{eq: Cradal Rabi}
Let $X$  and $Y$ be two Banach spaces, $U\subset X$ an open set of
$X$ and $I$ an open interval  of $\R$. We assume that  $0\in U$.
Denote by  $\varphi$ the elements of $U$  and  $\lambda$ the
elements of  $I$. Let $F: I\times U\rightarrow Y$  be a twice
continuously differentiable function such that
 \begin{enumerate}
\item
$ F(\lambda, 0)=0\quad \textrm{ for all }\quad \lambda \in I,$
 \item
 $\ker(D_\varphi F(\lambda_*, 0))=\R \varphi_*$  for some $\lambda_*\in I$  and $\varphi_*\in X\setminus \{0\}$,
 \item
 $\textrm{Codim Im}(D_\varphi F(\lambda_*, 0))=1,$
 \item
 $D_\lambda D_\varphi F(\lambda_*, 0)(\varphi_*) \notin \textrm{Im}(D_\varphi F(\lambda_*, 0)$.
 \end{enumerate}
 Then for any complement $Z$ of the subspace $\R \varphi_*$,  spanned by  $\varphi_*$, there exists a continuous curve
   $$(-\e, \e)\longrightarrow \R\times Z, \quad s\mapsto (\lambda(s), \varphi(s))$$  such that
 \begin{enumerate}
 \item
 $\lambda(0)=\lambda_*, \quad \varphi(0)=0,$
 \item
 $s(\varphi_*+\varphi (s))\in U,$
 \item
 $F(\lambda(s), s(\varphi_*+\varphi (s))=0$.
 \end{enumerate}
 Moreover, the set of solutions to the equation  $F(\lambda, u)=0$ in a neighborhood of  $(\lambda_*, 0)$
 is given by  the curve $\{ (\lambda, 0), \lambda\in \R\}$  and  $\{ s(\varphi_*+\varphi (s)), s\in(-\e, \e)\}$.
\end{Theorem}

\end{document}